\documentclass[12pt, reqno, twoside]{amsart}
\usepackage{amsmath, amsthm, amscd, amsfonts, amssymb, graphicx, color}
\usepackage{mathrsfs}
\usepackage[bookmarksnumbered, colorlinks, plainpages]{hyperref}
\hypersetup{colorlinks=true,linkcolor=red, anchorcolor=green, citecolor=cyan, urlcolor=red, filecolor=magenta, pdftoolbar=true}
\usepackage{enumitem}
\newcommand{\disp}{\displaystyle} 
\usepackage{enumitem}
\newtheorem{theorem}{Theorem}[section]
\newtheorem{lemma}[theorem]{Lemma}
\newtheorem{proposition}[theorem]{Proposition}
\newtheorem{corollary}[theorem]{Corollary}
\theoremstyle{definition}
\newtheorem{definition}[theorem]{Definition}
\newtheorem{example}[theorem]{Example}

\theoremstyle{remark}
\newtheorem{remark}[theorem]{Remark}
\numberwithin{equation}{section}

\begin{document}

\setcounter{page}{1}

\title[Hardy and Rellich types inequalities on the Heisenberg group]{Generalised Hardy type and Rellich type inequalities on the Heisenberg group}

\author[A. Abolarinwa]{Abimbola Abolarinwa}

\author[M. Ruzhansky]{Michael Ruzhansky}



\address{Abimbola Abolarinwa: Department of Mathematics,
University of Lagos, Akoka, Lagos State,
Nigeria\\ 
and}
\address{Department of Mathematics, Logic and Discrete Mathematics, 
Ghent University,  Ghent,  Belgium.}
\email{\textcolor[rgb]{0.00,0.00,0.84}{a.abolarinwa1@gmail.com, aabolarinwa@unilag.edu.ng}}

\address{Michael Ruzhansky:  Department of Mathematics, Logic and Discrete Mathematics, 
Ghent University,  Ghent,  Belgium \\
and}
\address{School of Mathematical Sciences, Queen Mary University, London, United Kingdom.}
\email{\textcolor[rgb]{0.00,0.00,0.84}{Michael.Ruzhansky@ugent.be}}
 

\begin{abstract}
This paper is primarily devoted to a class of interpolation inequalities of Hardy and Rellich types on the Heisenberg group $\mathbb{H}^n$.   Consequently, several weighted Hardy type, Heisenberg-Pauli-Weyl uncertainty principle and Hardy-Rellich type inqualities are established on $\mathbb{H}^n$. Moreover, new weighted Sobolev type embeddings are derived. Finally,  an integral inequality for vector fields in a domain of the Heisenberg group is obtained, leading to several specific weighted Hardy type inequalities  by making careful choices of vector fields. 
\end{abstract}
\maketitle
              
 
\section{Introduction and Preliminaries}
\subsection{Introduction}
Dated back to the works of G. H. Hardy (1877--1947), F. Rellich (1906--1955) and W. Heisenberg (1901--1976), Hardy type, Hardy-Rellich and uncertainty principle inequalities have found numerous applications in Mathematical analysis, theory of partial differential equations, physics, quantum mechanics,  etc. See the books and monographs \cite{HLP,OK,Re,RS} for details,  and also \cite{BV,Amb, Dev,GrA} and the references therein. Owing to their usefulness, these classes of inequalities have undergone series of extensions and refinements in the past decades. The interest of this present paper is on a class of interpolation inequalities of Caffarelli-Kohn-Nirenberg (CKN) type on the Heisenberg group. Explicit derivations of various forms of Hardy, Hardy-Rellich   and uncertainty principle type inequalities follow consequently.

Let $\mathbb{R}^n$ be the $n$-dimensional Euclidean space  and $f \in C^\infty_0(\mathbb{R}^n)$ be a smooth function. The CKN inequality \cite{CKN} states that there exists a constant $C_K>0$ such that 
\begin{align}\label{e1}
\||x|^cf\|_{L^r(\mathbb{R}^n)} \le C_K\||x|^a\nabla f\|^\delta_{L^p(\mathbb{R}^n)} \||x|^bf\|^{(1-\delta)}_{{L^q(\mathbb{R}^n)}},
\end{align}
where $|x|$ is the Euclidean norm and $r>0$, $p,q\ge 1$, $a,b,c \in \mathbb{R}$ and $\delta \in (0,1]$ satisfy
$$\frac{1}{r}+\frac{c}{n},\ \frac{1}{p}+\frac{a}{n},\ \frac{1}{q}+\frac{b}{n} >0, \ c=\delta\sigma+(1-\delta)b$$
if and only if 
\begin{align*}
\frac{1}{r}+\frac{c}{n}=\delta\left(\frac{1}{p}+\frac{a-1}{n}\right)+(1-\delta)\left(\frac{1}{q}+\frac{b}{n}\right)
\end{align*}
with\ 
$0\le \delta-\sigma \ \ \text{if}\ \delta>0$ and 
$a-\sigma\le 1 \ \ \text{if} \ \frac{1}{r}+\frac{c}{n} = \frac{1}{p}+\frac{a-1}{n}$.

If $a=0$, $\delta=1$ and $r=p$ are chosen, \eqref{e1} is the Hardy inequality, while it is the $L^2$-Sobolev inequality if $a=0$, $\delta=1$, $p=2$ and $r=2n/(n-2)$.  Other interesting and useful inequalities contained in \eqref{e1} are Gagliardo-Nirenberg, Hardy-Sobolev, Heisenberg-Pauli-Weyl (HPW), Nash inequalities and so on.
The proof of \eqref{e1} is more involved and has been generalised to some other settings such as Riemannian manifolds, homogeneous groups, hyperbolic spaces, geodesic spheres, see \cite{AbA, ARY,ARY3,AX,Amb2,KO1,ORS,ST,WN,WL} and the references therein for details. By the application of integration by parts and elementary algebraic inequality, D. G. Costa obtained new and short proofs for the $L^2$-version of CKN inequalities with sharp constants. In \cite{Co1}, he proved (Hardy type) that
\begin{align}\label{e2}
C_H\int_{\mathbb{R}^n}\frac{|f|^2}{|x|^{a+b+1}} dx \le \left(\int_{\mathbb{R}^n}\frac{|\nabla f|^2}{|x|^{2a}}\right)^{\frac{1}{2}}\left(\int_{\mathbb{R}^n}\frac{|f|^2}{|x|^{2b}}\right)^{\frac{1}{2}}
\end{align}
for $f \in C^\infty_0(\mathbb{R}^n\setminus\{0\})$, $a,b \in \mathbb{R}$ and $C_H:=|n-(a+b+1)|/2$ is sharp, and $\nabla$ is the Euclidean gradient operator. The possible values of $a$ and $b$ which are required in determining the exact values of the optimal constants and the functions that achieve them are discussed in \cite{CC}.  While in \cite{Co2}, he proved (Hardy-Rellich type) that
\begin{align}\label{e3}
C_R\int_{\mathbb{R}^n}\frac{|\nabla f|^2}{|x|^{a+b+1}} dx \le \left(\int_{\mathbb{R}^n}\frac{|\Delta f|^2}{|x|^{2a}}\right)^{\frac{1}{2}}\left(\int_{\mathbb{R}^n}\frac{|\nabla f|^2}{|x|^{2b}}\right)^{\frac{1}{2}}
\end{align}
for $a+b+3 \le n$ and $C_R:=|n+a+b-1|/2$ is sharp and $\Delta$ is the Euclidean Laplacian.

Note that \eqref{e2} and \eqref{e3} contain as special cases various forms of inequalities (Hardy, HPW, Rellich types). This can be seen by choosing different values for $a$ and $b$. Later on Di, Jang, Shen and Jin \cite{DJSJ} extended \eqref{e2} and \eqref{e3} to the general $L^p$ case ($1<p<\infty$) by a more direct method and they obtained
\begin{align}\label{e4}
C_{H,p}\int_{\mathbb{R}^n}\frac{|f|^p}{|x|^{a+b+1}} dx \le \left(\int_{\mathbb{R}^n}\frac{|\nabla f|^p}{|x|^{ap}}\right)^{\frac{1}{p}}\left(\int_{\mathbb{R}^n}\frac{|f|^p}{|x|^{bq}}\right)^{\frac{1}{q}},
\end{align}
where $1<p,q<\infty$, $1/p+1/q=1$ and the constant $C_{H,p}=|n-(a+b+1)|/p$ is sharp, and
\begin{align}\label{e5}
C_{R,p}\int_{\mathbb{R}^n}\frac{|\nabla f|^p}{|x|^{a+b+1}} dx \le \left(\int_{\mathbb{R}^n}\frac{|\Delta_p f|^p}{|x|^{ap}}\right)^{\frac{1}{p}}\left(\int_{\mathbb{R}^n}\frac{|\nabla f|^p}{|x|^{bq}}\right)^{\frac{1}{q}}
\end{align}
with $\frac{p-n}{p-1}\le \beta:=a+b+1$, $C_{R,p}:=(n+\beta(p-1)-p)/p$ and $\Delta_p$ is the usual $p$-Laplacian.

There have been few literature extending these results into domains where Euclidean distance is not available. For example, Abolarinwa, Rauf and Yin \cite{ARY3} recently obtained  versions of \eqref{e4} and \eqref{e5} on the unit $n$-sphere $\mathbb{S}^n$, using geodesic distance $r_q$ from a fixed point $q \in \mathbb{S}^n$, replacing $|x|$ with $|\tan r_q|$. Their inequalities have sharp constants with remainder terms. Wei and Li \cite{WL} extended \eqref{e4} into Riemannian manifolds using Hessian comparison theorem via construction of certain vector fields involving radial vector fields. 

Let $L^2_\alpha(\mathbb{R}^n)$, $D^{1,2}_\alpha(\mathbb{R}^n)$, $H^1_{a,b}(\mathbb{R}^n)$ and $H^2_{a,b}(\mathbb{R}^n)$ be the completions of $C^\infty_0(\mathbb{R}^n)$ respectively with the weighted Sobolev type norms\\
$\disp \|f\|_{L^2_\alpha(\mathbb{R}^n)}:= \left(\int_{\mathbb{R}^n }\frac{|f|^2}{|x|^{2\alpha}} dx \right)^{\frac{1}{2}}$,  
\hspace{1cm} $\disp \|f\|_{H^1_{a,b}(\mathbb{R}^n)}:= \left(\int_{\mathbb{R}^n }\left(\frac{|\nabla f|^2}{|x|^{2a}} + \frac{|f|^2}{|x|^{2b}}\right) dx \right)^{\frac{1}{2}}$, 
$\disp \|f\|_{D^{1,2}_\alpha(\mathbb{R}^n)}:= \left(\int_{\mathbb{R}^n }\frac{|\nabla f|^2}{|x|^{2\alpha}} dx \right)^{\frac{1}{2}}$,
and $\disp \|f\|_{H^2_{a,b}(\mathbb{R}^n)}:= \left(\int_{\mathbb{R}^n }\left(\frac{|\Delta f|^2}{|x|^{2a}} + \frac{|\nabla f|^2}{|x|^{2b}}\right) dx \right)^{\frac{1}{2}}$.
Then another consequences of the above equations \eqref{e2} and \eqref{e3} are the following continuous weighted Sobolev type embeddings \cite{Co1,Co2}:\\
$\disp H^1_{a,b}(\mathbb{R}^n) \subset L^2_{(a+b+a)/2}(\mathbb{R}^n)$, $\disp D^{2,2}_\alpha(\mathbb{R}^n) \subset D^{1,2}_{\alpha+1}(\mathbb{R}^n)$ and $\disp H^2_{a,b}(\mathbb{R}^n) \subset D^{1,2}_{(a+b+1)/2}(\mathbb{R}^n)$.

The main aim  of this paper is to extend \eqref{e2} and \eqref{e3} (resp. \eqref{e4} and \eqref{e5}) to the general $L^p$-version ($1<p<\infty$) in the setting of the Heisenberg group $\mathbb{H}^n$ using a more direct approach.  Since $\mathbb{H}^n$ is a stratified nilpotent Lie group, some methods used in the Euclidean setting cannot be imported directly to derive interpolation inequalities. The approach adopted here involved horizontal (Heisenberg norm) distance function (see Section \ref{sec2}).  Moreover, new weighted continuous Sobolev type embeddings are constructed by using Costa'a approach (Section \ref{sec3}). 
In Section \ref{sec4}, we obtain an integral inequality for vector fields in a domain of the Heisenberg group based on the integrability of vector fields and horizontal divergence theorem.  Finally, some specific Hardy type inequalities are then rederived by making careful choices of vector fields. 
The Heisenberg group is an important model to study many significant problems such as existence, continuation, eigenvalue problems, etc., so, obtaining vital integral inequalities such as those discussed in this paper become essential on $\mathbb{H}^n$. 

\subsection{Basics of the Heisenberg group}
Heisenberg group arises from many fields such as quantum physics, differential geometry, harmonic analysis. From the mathematical point of view, the Heisenberg group is the simplest and most important model in the theory of noncommutative vector fields satisfying  H\"orman\-der's condition. Here below we give basics of the group but detailed description can be found in \cite{FS,RS}.

Let $\xi = (x,y,t)\in \mathbb{R}^n\times \mathbb{R}^n\times \mathbb{R}$, $n\ge 1$. The set $\mathbb{H}^n=\mathbb{R}^{2n}\oplus\mathbb{R}$ equipped with the group law $\circ$,
\begin{align}\label{e6}
\xi\circ\xi_0=\Big(x+x_0,y+y_0, t+t_0+2\sum_{j=1}^n(x_{0j}\cdot y_j-x_j\cdot y_{0j})\Big),
\end{align}
where $x\cdot y$ is the usual Euclidean inner product in $\mathbb{R}^n$, is called the Heisenberg group of homogeneous dimension $Q=2n+2$. With the group law, $\mathbb{H}^n$ is a unimodular, connected and simply connected step two Lie group whose Haar measure is the Lebesgue measure, $d\xi=dxdydt$.
The distance between points $\xi$ and $\xi_0$ in $\mathbb{H}^n$ is defined by 
\begin{align}\label{e7}
d_{\mathbb{H}^n}(\xi,\xi_0)=|\xi_0^{-1}\circ\xi|_{\mathbb{H}^n},
\end{align}
where $\xi_0^{-1}$ is the inverse of $\xi_0$ given by $\xi_0^{-1}=-\xi_0$ with respect to the group law $\circ$. Let $z=(x,y)$, $\xi=(z,t)$ and $\xi_0=(z_0,t)$. By  \eqref{e6} and \eqref{e7}, one has
\begin{align}\label{e8}
d_{\mathbb{H}^n}(\xi,\xi_0)= \left([(x-x_0)^2+(y-y_0)^2]^2+[t-t_0-2(x\cdot y_0-x_0\cdot y)]^2 \right)^{\frac{1}{4}}.
\end{align}
The Heisenberg distance function (also called the Heisenberg norm) is therefore defined by
$$d:=d_{\mathbb{H}^n}(z,t)=(|z|^4+t^2)^{\frac{1}{4}}.$$
A basis for the Lie algebra of the left invariant vector fields on $\mathbb{H}^n$ is given by ($1\le j \le n$)
$$X_j=\frac{\partial}{\partial x_j}+2y_j\frac{\partial}{\partial t}, \ \ Y_j=\frac{\partial}{\partial y_j}-2x_j\frac{\partial}{\partial t}, \ \ 
 T= \frac{\partial}{\partial t}.$$
 By the Lie bracket definition one can easily check that 
 $$[X_j,X_k]=[Y_j,Y_k]=[X_j,T]=[Y_j,T]=0, \ j,k=1,2 \cdots n,$$
 and
 $$[X_j,Y_k]=-4T\delta_{jk},$$
which constitute Heisenberg canonical commutation relation of quantum mechanics for position and momentum, whence the name Heisenberg group.

A family of dilations is defined by $\delta_\lambda(z,t)=(\lambda z, \lambda^2t), \lambda>0$. It is also easy to check that $|\cdot|_{\mathbb{H}^n}, X_j, Y_j, (j=1,2 \cdots n)$ are homogeneous of degree $1$ with respect to dilations. 

The horizontal gradient, divergence and sub-Laplace operators on $\mathbb{H}^n$ are respectively defined by 
$$\nabla_{\mathbb{H}}f := (X_1f, \cdots, X_nf, Y_1f, \cdots , Y_nf), \ \ \text{div}_{\mathbb{H}} w:=\nabla_{\mathbb{H}} \cdot w \ \ \text{and} $$

$$\Delta_{\mathbb{H}}:=\sum_{j=1}^n\left(X_j^2+Y_j^2\right) = \Delta_z+4|z|^2\frac{\partial^2}{\partial t^2} + 4P\frac{\partial}{\partial t},$$
where $\Delta_z := \sum_{j=1}^n\left(\frac{\partial^2}{\partial x_j^2} +\frac{\partial^2}{\partial y_j^2}\right)$, $|z|^2= \sum_{j=1}^n\left(x_j^2+y_j^2\right)$ and $P:= \sum_{j=1}^n\left(y_j\frac{\partial}{\partial x_j}-x_j\frac{\partial}{\partial y_j}\right)$.
For $p>1$, the $p$-sub-Laplacian is also defined on $\mathbb{H}^n$ by 
$$\Delta_{\mathbb{H},p} f:= \text{div}_{\mathbb{H}}(|\nabla_{\mathbb{H}} f|^{p-2}\nabla_{\mathbb{H}} f).$$

The Heisenberg group belongs to a class of stratified Lie groups (or homogeneous Carnot groups). It was shown by Folland \cite{Fo} (see also \cite{FS}) that the sub-Laplacian on a general stratified Lie group has a unique fundamental solution $u_\varepsilon$ in distributional sense, which in the Heisenberg setting means
\begin{align*}
- \Delta_\mathbb{H} u_\varepsilon = \delta_0,
\end{align*}
where $\delta_0$ is the Dirac delta-distribution at the neutral element $0$ of $\mathbb{H}^n$.  The function (called $\Delta_\mathbb{H}$-gauge)
\begin{align*}
d(\xi):= \left\{
\begin{array}{ll}
\disp  u_\varepsilon(\xi)^{\frac{1}{2-Q}}, \ \  &\text{for} \ \xi \neq 0, \\ \ \\
\disp  0, \ & \text{for} \ \xi = 0,
\end{array}
\right.
\end{align*}
is a homogeneous quasi-norm on $\mathbb{H}^n$, that is, it is a continuous function smooth away from the origin.  The action of $\Delta_\mathbb{H}$ on $d$ is given as 
\begin{align*}
\Delta_\mathbb{H}d = (Q-1) \frac{|\nabla_\mathbb{H} d|^2}{d}  \ \ \text{in} \ \mathbb{H}^n\setminus\{0\}.
\end{align*}
We shall make use of the fact that Heisenberg is polarizable as a Carnot group \cite{RS},  such that  the $\infty$-sup-Laplacian of $d$ is 
\begin{align}\label{e9a}
\frac{1}{2} \langle \nabla_\mathbb{H} |\nabla_\mathbb{H} d|^2, \nabla_\mathbb{H} d \rangle =0  \ \ \text{in} \ \mathbb{H}^n\setminus\{0\}.
\end{align}

A simple computation therefore gives the following identities  in $\mathbb{H}^n\setminus\{0\}$ in explicit form:
\begin{align}\label{e9}
|\nabla_{\mathbb{H}} d| & =|z|d^{-1},
\end{align}
\begin{align}\label{e10}
\Delta_{\mathbb{H}} d &=  (Q-1)|z|^2d^{-3},
\end{align}
\begin{align}\label{e11}
\nabla_{\mathbb{H}}(|z|^{p-2})\nabla_{\mathbb{H}}d & =(p-2)|z|^pd^{-3}.
\end{align}
Using these identities one can easily compute
\begin{align}\label{e12}
\nabla_{\mathbb{H}}\left(\frac{|z|^{p-2}}{d^{p-2}}\right )\nabla_{\mathbb{H}}d =0.
\end{align}
Note that the standard Lebesgue measure on $\mathbb{R}^n$ is the Haar measure for the homogeneous group, so we write $d\xi=dzdt$, the Lebesgue measure on $\mathbb{H}^n$.

\section{Interpolation inequalities of Hardy and Rellich types}\label{sec2}
In this section we consider interpolation inequalities of Hardy and Rellich types.
\subsection{Sharp interpolation inequalities of Hardy type}
\begin{theorem}\label{thm21}
Let $\mathbb{H}^n$ be the Heisenberg group of homogeneous dimension $Q=2n+2$. Let $1<p<Q$ and $\lambda=a+b$, $a,b \in \mathbb{R}$ such that $Q\neq 1+\lambda$. Then we have
\begin{align}\label{e21}
\mathscr{C}_{a,b,p}^H \int_{\mathbb{H}^n} \frac{|z|^p}{d^p}\frac{|f|^p}{d^{\lambda+1}}d\xi  \le \left( \int_{\mathbb{H}^n} \frac{|\nabla_{\mathbb{H}}f|^p}{d^{ap}}d\xi\right)^{\frac{1}{p}}
\left(\int_{\mathbb{H}^n} \frac{|z|^p}{d^p}\frac{|f|^p}{d^{bq}}d\xi\right)^{\frac{1}{q}}
\end{align}
for  all $f \in D^{1,p}_a(\mathbb{H}^n)$,  with $D^{1,p}_a(\mathbb{H}^n)$ being the closure of  $C^\infty_0(\mathbb{H}^n\setminus\{0\})$ with respect to the norm
$$\|f\|_{D^{1,p}_a(\mathbb{H}^n)}:= \int_{\mathbb{H}^n} \frac{|\nabla_{\mathbb{H}}f|^p}{d^{ap}}d\xi,$$ 
where $\frac{1}{p}+\frac{1}{q}=1$ and $\mathscr{C}_{a,b,p}^H=|\frac{Q-1-\lambda}{p}|$ is sharp.
\end{theorem}

\proof
By a simple computation
\begin{align*}
\text{div}_{\mathbb{H}}\left(\frac{\nabla_{\mathbb{H}}d}{d^\lambda} \right) & = \frac{\Delta_{\mathbb{H}}d}{d^\lambda} - \frac{\lambda|\nabla_{\mathbb{H}}d|^2}{d^{\lambda+1}} \\
&= \frac{(Q-1)}{d^\lambda}\frac{|z|^2}{d^3} - \frac{\lambda}{d^{\lambda+1}}\frac{|z|^2}{d^2} = \frac{(Q-1-\lambda)}{d^{\lambda+1}}\frac{|z|^2}{d^2},
\end{align*}
where we have used \eqref{e9} and \eqref{e10}.  In the case $Q=\lambda+1$, the inequality is trivial, so we assume that  $Q\neq \lambda+1$. Applying the divergence theorem we have
\begin{align}\label{e22}
\int_{\mathbb{H}^n} \frac{|z|^p}{d^p}\frac{|f|^p}{d^{\lambda+1}}d\xi & = \frac{1}{Q-1-\lambda}\int_{\mathbb{H}^n} |f|^p\frac{|z|^{p-2}}{d^{p-2}} \nonumber \text{div}_{\mathbb{H}}\left(\frac{\nabla_{\mathbb{H}}d}{d^\lambda} \right) d\xi \\
& = \frac{-1}{Q-1-\lambda}\int_{\mathbb{H}^n} \nabla_{\mathbb{H}}\left(|f|^p\frac{|z|^{p-2}}{d^{p-2}}\right)\frac{\nabla_{\mathbb{H}}d}{d^\lambda} d\xi\\
&= \frac{-p}{Q-1-\lambda}\Re{e}\int_{\mathbb{H}^n}  |f|^{p-2} f\overline{\nabla_{\mathbb{H}}f}   \frac{\nabla_{\mathbb{H}}d}{d^\lambda} \frac{|z|^{p-2}}{d^{p-2}} d\xi, \nonumber
\end{align}
where we have applied \eqref{e12} into the second line to arrive at the last equality. On the other hand, application of Cauchy and H\"older inequalities gives
\begin{align*}
\frac{-p}{Q-1-\lambda}\Re{e}\int_{\mathbb{H}^n}  & |f|^{p-2}f\overline{\nabla_{\mathbb{H}}f}   \frac{\nabla_{\mathbb{H}}d}{d^\lambda} \frac{|z|^{p-2}}{d^{p-2}} d\xi \\ 
& \le \left|\frac{p}{Q-1-\lambda}\right|\int_{\mathbb{H}^n}  |f|^{p-1} |\nabla_{\mathbb{H}}f |  \frac{|\nabla_{\mathbb{H}}d|}{d^\lambda} \frac{|z|^{p-2}}{d^{p-2}} d\xi\\
& = \left|\frac{p}{Q-1-\lambda}\right|\int_{\mathbb{H}^n}  |f|^{p-1}   \frac{|\nabla_{\mathbb{H}}f |}{d^\lambda} \frac{|z|^{p-1}}{d^{p-1}} d\xi\\
&\le \left|\frac{p}{Q-1-\lambda}\right| \left(\int_{\mathbb{H}^n} \frac{|\nabla_{\mathbb{H}}f|^p}{d^{ap}}d\xi \right)^{\frac{1}{p}}\left(\int_{\mathbb{H}^n} \frac{|z|^p}{d^p}\frac{|f|^p}{d^{(\lambda-a)\frac{p}{p-1}}}d\xi\right)^{\frac{p-1}{p}}.
\end{align*}
Inserting the last inequality into \eqref{e22} leads to the required inequality \eqref{e21}.

To show the sharpness of the constant, we examine the equality condition in the H\"older inequality  (see \cite{ORS}). Define
\begin{align*}
g(x):=
\left\{
\begin{array}{ll}
e^{-\frac{C}{\beta}d^\beta} :\hspace{1cm} \beta=a-\frac{b}{p-1}+1 \neq 0 \\
d^{-C} \ \ : \hspace{1cm} \beta=a-\frac{b}{p-1}+1=0,
\end{array}
\right.
\end{align*}
where $C:=|\frac{Q-(\lambda+1)}{p}|$ and $Q\neq \lambda+1$. Then it can be checked that 
\begin{align*}
\Big|\frac{p}{Q-1-\lambda}\Big|^p \frac{|\nabla_{\mathbb{H}}g|^p}{d^{ap}} = \frac{|z|^p}{d^p} \frac{|g|^p}{d^{b\frac{p}{p-1}}},
\end{align*}
which satisfies the equality condition in the H\"older inequality. This shows that the constant $\mathscr{C}_{a,b,p}^H=|\frac{Q-1-\lambda}{p}|$ is sharp.
\qed

\begin{corollary}\label{cor22}
Let $\mathbb{H}^n$ be the Heisenberg group of homogeneous dimension $Q=2n+2$.  Let $1<p<Q$. Then for all $f \in D^{1,p}(\mathbb{H}^n)$  the following weighted Hardy inequality
\begin{align}\label{e23}
\left(\frac{Q-p}{p}\right)^p \int_{\mathbb{H}^n} \frac{|z|^p}{d^p}\frac{|f|^p}{d^p}d\xi  \le  \int_{\mathbb{H}^n} |\nabla_{\mathbb{H}}f|^p d\xi
\end{align}
holds with $(\frac{Q-p}{p})^p$ being the best constant.

Furthermore, the following HPW uncertainty principle inequality
\begin{align}\label{e24}
\frac{Q-p}{p} \int_{\mathbb{H}^n} \frac{|z|^p}{d^p}\frac{|f|^p}{d^p}d\xi \le \left(\int_{\mathbb{H}^n} |\nabla_{\mathbb{H}}f|^p d\xi \right)^{\frac{1}{p}} \left(\int_{\mathbb{H}^n}\frac{|z|^q}{d^q}\frac{|f|^q}{d^q} d\xi \right)^{\frac{1}{q}}
\end{align}
holds for $\frac{1}{p}+\frac{1}{q}=1$.  
\end{corollary}

\proof
Setting $a=0$ and $b=p-1$ in \eqref{e21} then $f\in D_0^{1,p}(\mathbb{H}^n)=D^{1,p}(\mathbb{H}^n)$ and \eqref{e21} yields \eqref{e23}.

Furthermore, by H\"older inequality (with conjugates $p$ and $q$) and \eqref{e23} we have
\begin{align*}
\int_{\mathbb{H}^n} \frac{|z|^p}{d^p}\frac{|f|^p}{d^p}d\xi &\le \left(\int_{\mathbb{H}^n} \frac{|z|^p}{d^p}\frac{|f|^p}{d^p}d\xi\right)^{\frac{1}{p}}\left(\int_{\mathbb{H}^n}\frac{|z|^q}{d^q}\frac{|f|^q}{d^q} d\xi \right)^{\frac{1}{q}}\\
&\le \frac{p}{Q-p}\left(\int_{\mathbb{H}^n} |\nabla_{\mathbb{H}}f|^p d\xi \right)^{\frac{1}{p}} \left(\int_{\mathbb{H}^n}\frac{|z|^q}{d^q}\frac{|f|^q}{d^q} d\xi \right)^{\frac{1}{q}},
\end{align*}
yielding \eqref{e24}.
\qed

In the Euclidean setting, the error term between left hand and right hand side of \eqref{e23} was investigated by Ioku, Ishiwata and Ozawa \cite{IIO}. See \cite{RS} for the case of general homogeneous group with radial operator and quasi-norms. Note also that \eqref{e23} has been derived by Niu, Zhang and Wang \cite{NZW} through the application of generalised Picone formula to the existence of subsolution of $p$-subLaplacian eigenvalue.  A similar approach to \cite{NZW} was adopted by Han and Niu \cite{HN} to obtain
\begin{align*}
\mathscr{D}\int_\Omega \frac{|z|^p}{d^p}\frac{|f|^p}{d^{2p}}d\xi \le \int_\Omega |\nabla_{\mathbb{H}}f|^p d\xi,
\end{align*}
where the constant $\mathscr{D}>0$, $\Omega = \mathbb{B}_R\setminus\{0\}$ and $\mathbb{B}_R=\{(z,t)\in \mathbb{H}^n : d(z,t)<R\}$.  

The weighted Hardy inequality \eqref{e23} also gives a special case of Han, Niu and Zhang \cite[Theorem 3.2]{HNZ} with sharp constant derived via a representation formula for functions. The case $p=2$ in \eqref{e23} (that is when $a=0$, $b=0$ and $p=2$ for $f \in D^{1,2}(\mathbb{H}^n)$) has been discussed by Garofalo and Lanconelli \cite{GL} in their Corollaries 2.1 and 2.2 using a different method which involves the use of the fundamental solution of $-\Delta_{\mathbb{H}}$. They derived an uncertainty principle
\begin{align*}
\left(\frac{Q-2}{2}\right)^2\left(\int_{\mathbb{H}^n}\frac{|z|^2}{d^2}|f|^2d\xi\right)^2 \le \left(\int_{\mathbb{H}^n}|z|^2|f|^2d\xi\right)\left(\int_{\mathbb{H}^n}|\nabla_{\mathbb{H}^n}f|^2d\xi\right)
\end{align*}
and remarked that equality holds when $(Q-2)/2$ is replaced by $Q/2$ if and only if $f(z,t)=A\exp(-\alpha d^2(z,t))$ ($A\in \mathbb{R}$, $\alpha>0$). 

Several uncertainty principle inequalities can be derived directly from \eqref{e21} by specializing $a$ and $b$. The next Corollary is given to demonstrate this.

\begin{corollary}\label{cor23}
Given the assumptions of Theorem \ref{thm21}.
\begin{enumerate}[label=(\roman*)]
\item When $a+b=0$ and $a=-p$, then for $Q>1$
\begin{align*}
\frac{Q-1}{p}\int_{\mathbb{H}^n} \frac{|z|^p}{d^p}\frac{|f|^p}{d^p}d\xi \le \left(\int_{\mathbb{H}^n} |\nabla_{\mathbb{H}}f|^p d^{p^2} d\xi \right)^{\frac{1}{p}} \left(\int_{\mathbb{H}^n}\frac{|z|^p}{d^p}\frac{|f|^p}{d^{pq}} d\xi \right)^{\frac{1}{q}}.
\end{align*}
\item When $a+b+1=0$ and $a=-p$, then 
\begin{align*}
\frac{Q}{p}\int_{\mathbb{H}^n} \frac{|z|^p}{d^p}|f|^pd\xi \le \left(\int_{\mathbb{H}^n} |\nabla_{\mathbb{H}}f|^p d^{p^2} d\xi \right)^{\frac{1}{p}} \left(\int_{\mathbb{H}^n}\frac{|z|^p}{d^p}\frac{|f|^p}{d^{p}} d\xi \right)^{\frac{1}{q}}.
\end{align*}
\item When $a+b+1=0$ and $a=0$, then 
\begin{align*}
\frac{Q}{p}\int_{\mathbb{H}^n} \frac{|z|^p}{d^p}|f|^pd\xi \le \left(\int_{\mathbb{H}^n} |\nabla_{\mathbb{H}}f|^p d\xi \right)^{\frac{1}{p}} \left(\int_{\mathbb{H}^n}\frac{|z|^p}{d^p}|f|^p d^q d\xi \right)^{\frac{1}{q}}.
\end{align*}
\item When $a+b+1=0$ and $a=-2$, then 
\begin{align*}
\frac{Q}{p}\int_{\mathbb{H}^n} \frac{|z|^p}{d^p}|f|^pd\xi \le \left(\int_{\mathbb{H}^n} |\nabla_{\mathbb{H}}f|^p d^{2p}d\xi \right)^{\frac{1}{p}} \left(\int_{\mathbb{H}^n}\frac{|z|^p}{d^p}\frac{|f|^p}{d^p} d\xi \right)^{\frac{1}{q}}.
\end{align*}
\item When $a+b+1=0$ and $a=1$, then 
\begin{align*}
\frac{Q}{p}\int_{\mathbb{H}^n} \frac{|z|^p}{d^p}|f|^pd\xi \le \left(\int_{\mathbb{H}^n} \frac{|\nabla_{\mathbb{H}}f|^p}{d^p} d\xi \right)^{\frac{1}{p}} \left(\int_{\mathbb{H}^n}\frac{|z|^p}{d^p}|f|^p d^{2q} d\xi \right)^{\frac{1}{q}}.
\end{align*}
\item When $a+b+1=ap$ and $Q\neq ap$, then 
\begin{align*}
\frac{|Q-ap|}{p}\int_{\mathbb{H}^n} \frac{|z|^p}{d^p}\frac{|f|^p}{d^{ap}}d\xi \le \left(\int_{\mathbb{H}^n} \frac{|\nabla_{\mathbb{H}}f|^p}{d^{ap}} d\xi \right)^{\frac{1}{p}} \left(\int_{\mathbb{H}^n}\frac{|z|^p}{d^p}\frac{|f|^p}{d^{ap}} d^q d\xi \right)^{\frac{1}{q}}.
\end{align*}
\end{enumerate}
\end{corollary}

\subsection{Rellich type interpolation inequalities}
\begin{theorem}\label{thm24}
Let $\mathbb{H}^n$ be the Heisenberg group of homogeneous dimension $Q=2n+2$. Let $1<p<Q$ and $\lambda=a+b$  be such that for $a,b \in \mathbb{R}$,  $\frac{p-Q}{p-1}\le \lambda+1\le 0$. Then  for all $f \in C^\infty_0(\mathbb{H}^n\setminus\{0\})$ we have
\begin{align}\label{e25}
\mathscr{C}_{a,b,p}^R \int_{\mathbb{H}^n} \frac{|z|^p}{d^p}\frac{|\nabla_{\mathbb{H}}f|^p}{d^{\lambda+1}}d\xi  \le \left( \int_{\mathbb{H}^n} \frac{|\Delta_{\mathbb{H},p}f|^p}{d^{ap}}d\xi\right)^{\frac{1}{p}}
\left(\int_{\mathbb{H}^n} \frac{|z|^p}{d^p}\frac{|\nabla_{\mathbb{H}}f|^{\frac{p}{p-1}}}{d^{b\frac{p}{p-1}}}d\xi\right)^{\frac{p-1}{p}},
\end{align} 
where the constant is $\mathscr{C}_{a,b,p}^R:=\frac{Q-p+(\lambda+1)(p-1)}{p}$  and $\Delta_{\mathbb{H},p}$ is the $p$-sub-Laplacian defined by  $\Delta_{\mathbb{H},p} f:= \text{div}_{\mathbb{H}}(|\nabla_{\mathbb{H}} f|^{p-2}\nabla_{\mathbb{H}} f).$
\end{theorem}

\begin{remark}
It is not clear whether the constant $\mathscr{C}_{a,b,p}^R:=\frac{Q-p+(\lambda+1)(p-1)}{p}$  is optimal or not.
\end{remark}

\proof
Recall that a simple computation gives
\begin{align*}
\text{div}_{\mathbb{H}}\left(\frac{\nabla_{\mathbb{H}}d}{d^\lambda} \right) = \frac{(Q-1-\lambda)}{d^{\lambda+1}}\frac{|z|^2}{d^2}.
\end{align*}
Then by the divergence theorem and \eqref{e12} we get
\begin{align}\label{e26}
\int_{\mathbb{H}^n} \frac{|z|^p}{d^p}\frac{|\nabla_{\mathbb{H}}f|^p}{d^{\lambda+1}}d\xi = \frac{-p}{Q-1-\lambda}  \int_{\mathbb{H}^n}  \frac{|z|^{p-2}}{d^{p-2}} |\nabla_{\mathbb{H}}f|^{p-2} \frac{ \frac{1}{2} \nabla_{\mathbb{H}}(|\nabla_\mathbb{H}f |^2) \nabla_{\mathbb{H}}d}{d^\lambda}  d\xi.
\end{align}
On the other hand, also by the divergence theorem and \eqref{e12} we have
\begin{align*}
\Re{e} \int_{\mathbb{H}^n} & \Delta_{\mathbb{H},p}f \frac{\overline{\nabla_{\mathbb{H}}f} \nabla_{\mathbb{H}}d}{d^\lambda} \frac{|z|^{p-2}}{d^{p-2}} d\xi \\
 & = - \Re{e}\int_{\mathbb{H}^n}   |\nabla_{\mathbb{H}}f|^{p-2} \nabla_{\mathbb{H}}f \nabla_{\mathbb{H}} \left(\frac{\overline{\nabla_{\mathbb{H}}f} \nabla_{\mathbb{H}}d}{d^\lambda} \frac{|z|^{p-2}}{d^{p-2}} \right) d\xi\\
&=- \Re{e} \int_{\mathbb{H}^n}   |\nabla_{\mathbb{H}}f|^{p-2} \nabla_{\mathbb{H}}f \nabla_{\mathbb{H}} \left(\frac{\overline{\nabla_{\mathbb{H}}f} \nabla_{\mathbb{H}}d}{d^\lambda} \right) \frac{|z|^{p-2}}{d^{p-2}}  d\xi\\
&=-\Re{e}\int_{\mathbb{H}^n}   |\nabla_{\mathbb{H}}f|^{p-2} \Bigg( \nabla_{\mathbb{H}}f \frac{ \nabla_{\mathbb{H}}(\overline{\nabla_{\mathbb{H}}f}) \nabla_{\mathbb{H}}d}{d^\lambda} + \frac{\nabla_{\mathbb{H}}f\overline{\nabla_{\mathbb{H}}f}\nabla_{\mathbb{H}}(\nabla_{\mathbb{H}}d)}{d^\lambda} \\ 
&\hspace{4.5cm} - \frac{\lambda |\nabla_{\mathbb{H}}f|^2 |\nabla_{\mathbb{H}}d|^2}{d^{\lambda+1}}
\Bigg)\frac{|z|^{p-2}}{d^{p-2}}  d\xi\\
&=-\int_{\mathbb{H}^n}   |\nabla_{\mathbb{H}}f|^{p-2} \left( \frac{ \frac{1}{2} \nabla_{\mathbb{H}}(|\nabla_\mathbb{H}f |^2) \nabla_{\mathbb{H}}d}{d^\lambda} - \frac{\lambda |\nabla_{\mathbb{H}}f|^2}{d^{\lambda+1}}\frac{|z|^2}{d^2}\right)\frac{|z|^{p-2}}{d^{p-2}}  d\xi,
\end{align*}
where the middle term in the third equality vanishes in the sense of \eqref{e9a}. That is,
$$ \nabla_\mathbb{H}(\nabla_\mathbb{H}d ) \frac{|z|^{p-2}}{d^{p-2}} =   \nabla_\mathbb{H}(\nabla_\mathbb{H}d ) | \nabla_\mathbb{H}d|^{p-2} = \frac{1}{2} \langle \nabla_\mathbb{H} |\nabla_\mathbb{H} d|^2, \nabla_\mathbb{H} d \rangle | \nabla_\mathbb{H}d|^{p-4} =0.$$ 
 Therefore, 
\begin{align}\label{e27}
\int_{\mathbb{H}^n}  \frac{|z|^{p-2}}{d^{p-2}} |\nabla_{\mathbb{H}}f|^{p-2}  \frac{\frac{1}{2} \nabla_{\mathbb{H}}(|\nabla_{\mathbb{H}}f|^2) \nabla_{\mathbb{H}}d}{d^\lambda}  d\xi =& \lambda \int_{\mathbb{H}^n} \frac{|z|^p}{d^p} \frac{|\nabla_{\mathbb{H}}f|^p}{d^{\lambda+1}} d\xi \nonumber\\
& -\Re{e} \int_{\mathbb{H}^n}  \Delta_{\mathbb{H},p}f \frac{\overline{\nabla_{\mathbb{H}}f} \nabla_{\mathbb{H}}d}{d^\lambda} \frac{|z|^{p-2}}{d^{p-2}} d\xi.
\end{align}
Combining \eqref{e26} and \eqref{e27} yields
\begin{align*}
\int_{\mathbb{H}^n} \frac{|z|^p}{d^p}\frac{|\nabla_{\mathbb{H}}f|^p}{d^{\lambda+1}}d\xi = \frac{-p}{Q-1-\lambda} \Big(\lambda \int_{\mathbb{H}^n} \frac{|z|^p}{d^p} \frac{|\nabla_{\mathbb{H}}f|^p}{d^{\lambda+1}} d\xi -\Re{e} \int_{\mathbb{H}^n}  \Delta_{\mathbb{H},p}f \frac{\overline{\nabla_{\mathbb{H}}f} \nabla_{\mathbb{H}}d}{d^\lambda} \frac{|z|^{p-2}}{d^{p-2}} d\xi\Big)
\end{align*}
which after rearrangement implies
\begin{align}\label{e28}
\frac{Q-1+\lambda(p-1)}{p} \int_{\mathbb{H}^n} \frac{|z|^p}{d^p}\frac{|\nabla_{\mathbb{H}}f|^p}{d^{\lambda+1}}d\xi =\Re{e} \int_{\mathbb{H}^n}  \Delta_{\mathbb{H},p}f \frac{\overline{\nabla_{\mathbb{H}}f} \nabla_{\mathbb{H}}d}{d^\lambda} \frac{|z|^{p-2}}{d^{p-2}} d\xi.
\end{align}
By the condition $\frac{p-Q}{p-1}\le \lambda+1=a+b+1\le 0$, we note that\\ 
$\mathscr{C}_{a,b,p}^R:= \frac{Q-1+\lambda(p-1)}{p} = \frac{Q-p+(\lambda+1)(p-1)}{p}\ge 0.$

Now by the  H\"older inequality and \eqref{e9} we get
\begin{align*}
\Re{e}\int_{\mathbb{H}^n}  \Delta_{\mathbb{H},p}f \frac{\overline{\nabla_{\mathbb{H}}f}\nabla_{\mathbb{H}}d}{d^\lambda} \frac{|z|^{p-2}}{d^{p-2}} d\xi & \le \int_{\mathbb{H}^n} | \Delta_{\mathbb{H},p}f| \frac{|\nabla_{\mathbb{H}}f|}{d^\lambda} \frac{|z|^{p-1}}{d^{p-1}} d\xi\\
& \le \left( \int_{\mathbb{H}^n} \frac{|\Delta_{\mathbb{H},p}f|^p}{d^{ap}}d\xi\right)^{\frac{1}{p}}
\left(\int_{\mathbb{H}^n} \frac{|z|^p}{d^p}\frac{|\nabla_{\mathbb{H}}f|^q}{d^{(\lambda-a)q }}d\xi\right)^{\frac{1}{q}},
\end{align*}
with $\frac{1}{p}+\frac{1}{q}=1$.  Inserting the last inequality into \eqref{e28} yields the required inequality.
\qed

\begin{remark}\label{rem26}
We give several special cases of  \eqref{e27} by choosing specific values of $a$ and $b$, especially for $a+b+1=0$ ($1<p<Q$ and $\frac{1}{p}+ \frac{1}{q}=1$):
\begin{enumerate}[label=(\roman*)]
\item When $a=0$ and $b=-1$, 
\begin{align*}
\frac{Q-p}{p}\int_{\mathbb{H}^n} \frac{|z|^p}{d^p} |\nabla_{\mathbb{H}} f|^p d\xi \le \left(\int_{\mathbb{H}^n} |\Delta_{\mathbb{H},p}f|^p d\xi \right)^{\frac{1}{p}} \left(\int_{\mathbb{H}^n}\frac{|z|^p}{d^p}|\nabla_{\mathbb{H}} f|^q d^q d\xi \right)^{\frac{1}{q}}.
\end{align*}
\item When $a=-1$ and $b=0$,
\begin{align*}
\frac{Q-p}{p}\int_{\mathbb{H}^n} \frac{|z|^p}{d^p} |\nabla_{\mathbb{H}} f|^p d\xi \le \left(\int_{\mathbb{H}^n} |\Delta_{\mathbb{H},p}f|^p d^p d\xi \right)^{\frac{1}{p}} \left(\int_{\mathbb{H}^n}\frac{|z|^p}{d^p}|\nabla_{\mathbb{H}} f|^q d\xi \right)^{\frac{1}{q}}.
\end{align*}
A special case of this is when $p=2$:
\begin{align*}
\left(\frac{Q-2}{2}\right)^2\int_{\mathbb{H}^n} \frac{|z|^2}{d^2} |\nabla_{\mathbb{H}} f|^2 d\xi \le \int_{\mathbb{H}^n} |\Delta_{\mathbb{H}}f|^2 d^2 d\xi.
\end{align*}
\item When $a=1$ and $b=-2$,
\begin{align*}
\frac{Q-p}{p}\int_{\mathbb{H}^n} \frac{|z|^p}{d^p} |\nabla_{\mathbb{H}} f|^p d\xi \le \left(\int_{\mathbb{H}^n} \frac{|\Delta_{\mathbb{H},p}f|^p}{d^p} d\xi \right)^{\frac{1}{p}} \left(\int_{\mathbb{H}^n}\frac{|z|^p}{d^p}|\nabla_{\mathbb{H}} f|^q d^{2q}d\xi \right)^{\frac{1}{q}}.
\end{align*}
\item When $a=-2$ and $b=1$,
\begin{align*}
\frac{Q-p}{p}\int_{\mathbb{H}^n} \frac{|z|^p}{d^p} |\nabla_{\mathbb{H}} f|^p d\xi \le \left(\int_{\mathbb{H}^n} |\Delta_{\mathbb{H},p}f|^p d^{2p} d\xi \right)^{\frac{1}{p}} \left(\int_{\mathbb{H}^n}\frac{|z|^p}{d^p}\frac{|\nabla_{\mathbb{H}} f|^q}{d^q} d\xi \right)^{\frac{1}{q}}.
\end{align*}
\end{enumerate}
\end{remark}

\section{Weighted Sobolev  type embedding theorems}\label{sec3}
As a consequence of the horizontal weights in Theorem \ref{thm21} and Theorem \ref{thm24} we are prompted to define the following weighted Sobolev type spaces on the Heisenberg group $\mathbb{H}^n$ for $p=2$. We then establish several horizontal embeddings. Since we are now in the regime of $p=2$ we shall take the advantage to extend Costa's approach \cite{Co1,Co2} to our setting.

\begin{definition}{\bf (Weighted Sobolev type spaces)}
The following weighted Sobolev type spaces are defined  with respect to the following norms as the completion of $C^\infty_0(\mathbb{H}^n\setminus\{0\})$.
\end{definition}

\begin{enumerate}[label=(\alph*)]
\item $L^2_\alpha(\mathbb{H}^n,\frac{|z|^2}{d^2})$ denotes the completion of $C^\infty_0(\mathbb{H}^n\setminus\{0\})$ with the respect to the norm
$$\|f\|_{L^2_\alpha(\mathbb{H}^n,\frac{|z|^2}{d^2})}:=\left(\int_{\mathbb{H}^n} \frac{|z|^2}{d^2}\frac{|f|^2}{d^{2\alpha}}d\xi \right)^{\frac{1}{2}}.$$

\item $D^{1,2}_\alpha(\mathbb{H}^n)$ denotes the completion of $C^\infty_0(\mathbb{H}^n\setminus\{0\})$ with the respect to the norm
$$\|f\|_{D^{1,2}_\alpha(\mathbb{H}^n)}:=\left(\int_{\mathbb{H}^n} \frac{|\nabla_{\mathbb{H}}f|^2}{d^{2\alpha}}d\xi \right)^{\frac{1}{2}}.$$

\item $D^{1,2}_\alpha(\mathbb{H}^n,\frac{|z|^2}{d^2})$ denotes the completion of $C^\infty_0(\mathbb{H}^n\setminus\{0\})$ with the respect to the norm
$$\|f\|_{D^{1,2}_\alpha(\mathbb{H}^n,\frac{|z|^2}{d^2})}:=\left(\int_{\mathbb{H}^n} \frac{|z|^2}{d^2}\frac{|\nabla_{\mathbb{H}}f|^2}{d^{2\alpha}}d\xi \right)^{\frac{1}{2}}.$$

\item $D^{2,2}_\alpha(\mathbb{H}^n)$ denotes the completion of $C^\infty_0(\mathbb{H}^n\setminus\{0\})$ with the respect to the norm
$$\|f\|_{D^{2,2}_\alpha(\mathbb{H}^n)}:=\left(\int_{\mathbb{H}^n} \frac{|\Delta_{\mathbb{H}}f|^2}{d^{2\alpha}}d\xi \right)^{\frac{1}{2}}.$$

\item $H^1_{a,b}(\mathbb{H}^n,\frac{|z|^2}{d^2})$ denotes the completion of $C^\infty_0(\mathbb{H}^n\setminus\{0\})$ with the respect to the norm
$$\|f\|_{H^1_{a,b}(\mathbb{H}^n,\frac{|z|^2}{d^2})}:=\left(\int_{\mathbb{H}^n} \frac{|\nabla_{\mathbb{H}}f|^2}{d^{2a}}d\xi
+\int_{\mathbb{H}^n} \frac{|z|^2}{d^2}\frac{|f|^2}{d^{2b}}d\xi \right)^{\frac{1}{2}}.$$

\item $H^2_{a,b}(\mathbb{H}^n,\frac{|z|^2}{d^2})$ denotes the completion of $C^\infty_0(\mathbb{H}^n\setminus\{0\})$ with the respect to the norm
$$\|f\|_{H^2_{a,b}(\mathbb{H}^n,\frac{|z|^2}{d^2})}:=\left(\int_{\mathbb{H}^n} \frac{|\Delta_{\mathbb{H},p}f|^2}{d^{2a}}d\xi
+\int_{\mathbb{H}^n} \frac{|z|^2}{d^2}\frac{|\nabla_{\mathbb{H}}f|^2}{d^{2b}}d\xi \right)^{\frac{1}{2}}.$$
\end{enumerate}

\begin{theorem}\label{thm32}
{\bf (Sobolev type embedding with horizontal weights)}
Let $\mathbb{H}^n$ be the Heisenberg group of homogeneous dimension $Q=2n+2$ and let $a,b \in \mathbb{R}$.
\begin{enumerate}[label=(\alph*)]
\item For $Q\neq a+b+1$,  the following continuous embedding holds
\begin{enumerate}[label=(\roman*)]
\item $H^1_{a,b}(\mathbb{H}^n,\frac{|z|^2}{d^2}) \subset L^2_{\frac{a+b+1}{2}}(\mathbb{H}^n,\frac{|z|^2}{d^2})$.\\
Moreover,
\item $H^1_{b,a}(\mathbb{H}^n,\frac{|z|^2}{d^2}) \subset L^2_{\frac{a+b+1}{2}}(\mathbb{H}^n,\frac{|z|^2}{d^2})$
due to symmetry with respect to $a$ and $b$.
\end{enumerate}
\item For $Q\neq a+b-1$,  the following continuous embeddings hold
\begin{enumerate}[label=(\roman*)]
\item $H^2_{a,b}(\mathbb{H}^n,\frac{|z|^2}{d^2}) \subset D^{1,2}_{\frac{a+b+1}{2}}(\mathbb{H}^n,\frac{|z|^2}{d^2})$\\
and
\item $D^{2,2}_\alpha(\mathbb{H}^n) \subset D^{1,2}_{\alpha+1}(\mathbb{H}^n,\frac{|z|^2}{d^2})$
for $\alpha \le \frac{Q}{2}-2$ and $\alpha \neq \frac{Q}{2}$.
\end{enumerate}
\end{enumerate}
\end{theorem}

\subsection{Proof of Part (a) in Theorem \ref{thm32} and some consequences}
To prove part (a) Theorem \ref{thm32}, we start with the case $p=2$ of Theorem \ref{thm21}, that is, the inequality stated in the next Lemma.

\begin{lemma}\label{lem33}
For $Q\neq a+b+1$, $a,b \in \mathbb{R}$,
\begin{align}\label{e31}
\mathscr{D}_{a,b} \int_{\mathbb{H}^n} \frac{|z|^2}{d^2}\frac{|f|^2}{d^{a+b+1}}d\xi  \le \left( \int_{\mathbb{H}^n} \frac{|\nabla_{\mathbb{H}}f|^2}{d^{2a}}d\xi\right)^{\frac{1}{2}}
\left(\int_{\mathbb{H}^n} \frac{|z|^2}{d^2}\frac{|f|^2}{d^{2b}}d\xi\right)^{\frac{1}{2}}
\end{align}
with sharp constant $\mathscr{D}_{a,b}:=\frac{|Q-(a+b+1)|}{2}$.
\end{lemma}

\subsection*{Proof of Lemma \ref{lem33}} We discuss two approaches to this.\\
{\bf First Approach:} As stated above, setting $p=2$ in Theorem \eqref{thm21} leads to \eqref{e31}.\\
{\bf Second Approach:} - It follows from Costa's approach \cite{Co1}. For all $f \in C^\infty_0(\mathbb{H}^n\setminus\{0\})$ and $a,b, s \in \mathbb{R}$, we have
\begin{align}\label{e32}
\int_{\mathbb{H}^n} \left|\frac{\nabla_{\mathbb{H}}f}{d^a}+s\frac{f}{d^b}\nabla_{\mathbb{H}^n} d \right|^2 d\xi \ge 0
\end{align}
implying that
\begin{align}\label{e33}
\int_{\mathbb{H}^n} \frac{|\nabla_{\mathbb{H}} f|^2}{d^{2a}}d\xi
+ s^2 \int_{\mathbb{H}^n} \frac{|z|^2}{d^2}\frac{|f|^2}{d^{2b}}d\xi
+2s \Re{e}\int_{\mathbb{H}^n} \frac{f}{d^{a+b}} \overline{\nabla_{\mathbb{H}} f} \nabla_{\mathbb{H}} d d\xi \ge 0.
\end{align}
By the divergence theorem
\begin{align*}
\Re{e}\int_{\mathbb{H}^n} & \frac{f}{d^{a+b}} \overline{\nabla_{\mathbb{H}} f} \nabla_{\mathbb{H}} d d\xi 
 = - \Re{e}\int_{\mathbb{H}^n} \overline{f} \text{div}_{\mathbb{H}} \left(\frac{f}{d^{a+b}} \nabla_{\mathbb{H}} d \right) d\xi\\
=& - \int_{\mathbb{H}^n} \frac{|f|^2}{d^{a+b}}\Delta_{\mathbb{H}} d d\xi - \Re{e}\int_{\mathbb{H}^n} \frac{f}{d^{a+b}}\overline {\nabla_{\mathbb{H}}f} \nabla_{\mathbb{H}} d d\xi +(a+b) \int_{\mathbb{H}^n} \frac{|f|^2}{d^{a+b+1}} |\nabla_{\mathbb{H}} d|^2 d\xi\\
=& -(Q-(a+b+1)) \int_{\mathbb{H}^n} \frac{|f|^2}{d^{a+b+1}}\frac{|z|^2}{d^2} d\xi - \Re{e}\int_{\mathbb{H}^n}  \frac{f}{d^{a+b}}\overline{ \nabla_{\mathbb{H}} f} \nabla_{\mathbb{H}} d d\xi,
\end{align*}
where we have used \eqref{e9} and \eqref{e10}.  Thus,
\begin{align*}
2 \Re{e}\int_{\mathbb{H}^n} & \frac{f}{d^{a+b}} \overline{\nabla_{\mathbb{H}} f} \nabla_{\mathbb{H}} d d\xi = - (Q-(a+b+1))\int_{\mathbb{H}^n} \frac{|f|^2}{d^{a+b+1}}\frac{|z|^2}{d^2} d\xi.
\end{align*}
Furthermore, denote by \\
$$A:= \int_{\mathbb{H}^n} \frac{|z|^2}{d^2}\frac{|f|^2}{d^{2b}}d\xi, \ \ B:=  \Re{e}\int_{\mathbb{H}^n}  \frac{f}{d^{a+b}} \overline{\nabla_{\mathbb{H}} f} \nabla_{\mathbb{H}} d d\xi , \ \ 
 \text{and} \ \ \ C:= \int_{\mathbb{H}^n} \frac{|\nabla_{\mathbb{H}} f|^2}{d^{2a}}d\xi. $$
 Then \eqref{e33} takes the form $As^2+Bs+C\ge 0$ which holds true if and only if $B^2-4AC\le 0$. Therefore 
\begin{align*}
\left((Q-(a+b+1))\int_{\mathbb{H}^n} \frac{|f|^2}{d^{a+b+1}}\frac{|z|^2}{d^2} d\xi \right)^2 \le 4 \int_{\mathbb{H}^n} \frac{|\nabla_{\mathbb{H}} f|^2}{d^{2a}}d\xi \int_{\mathbb{H}^n} \frac{|z|^2}{d^2}\frac{|f|^2}{d^{2b}}d\xi.
\end{align*}
The proof is complete.

\subsection*{Proof of Part (a) in Theorem \eqref{thm32}}
Since $Q\neq a+b+1$, using Young's inequality in \ref{e31} one obtains
\begin{align*}
\int_{\mathbb{H}^n} \frac{|z|^2}{d^2}\frac{|f|^2}{d^{a+b+1}}d\xi  & \le \frac{2}{|Q-(a+b+1)|}\left( \int_{\mathbb{H}^n} \frac{|\nabla_{\mathbb{H}}f|^2}{d^{2a}}d\xi\right)^{\frac{1}{2}}
\left(\int_{\mathbb{H}^n} \frac{|z|^2}{d^2}\frac{|f|^2}{d^{2b}}d\xi\right)^{\frac{1}{2}}\\
& \le \frac{1}{|Q-(a+b+1)|}\left(\int_{\mathbb{H}^n} \frac{|\nabla_{\mathbb{H}}f|^2}{d^{2a}}d\xi + \int_{\mathbb{H}^n} \frac{|z|^2}{d^2}\frac{|f|^2}{d^{2b}}d\xi\right)
\end{align*}
for all $f \in C^\infty_0(\mathbb{H}^n\setminus\{0\})$. This implies the embedding (i). The arbitrariness of the real parameters $a$ and $b$ and symmetry of \eqref{e31} with respect to $a$ and $b$ yields (ii).
\qed

As a consequence of Lemma \ref{lem33}, various forms of generalised inequalities can be obtained.
\begin{corollary}\label{cor34}
The following inequalities hold and the constants are sharp.
\begin{enumerate}[label=(\arabic*)]
\item For any $f \in H^1_{a,a+1}(\mathbb{H}^n,\frac{|z|^2}{d^2})$, it holds that
\begin{align*}
\left(\frac{Q}{2}-(a+1)\right)^2 \int_{\mathbb{H}^n} \frac{|z|^2}{d^2}\frac{|f|^2}{d^{2(a+1)}}d\xi  \le \int_{\mathbb{H}^n} \frac{|\nabla_{\mathbb{H}}f|^2}{d^{2a}}d\xi
\end{align*}
and $D^{1,2}_\alpha(\mathbb{H}^n) \subset L^1_{\alpha+1}(\mathbb{H}^n,\frac{|z|^2}{d^2})$.
\item For any $f \in H^1_{b+1,b}(\mathbb{H}^n,\frac{|z|^2}{d^2})$, it holds that
\begin{align*}
\left(\frac{Q}{2}-(b+1)\right)^2 \int_{\mathbb{H}^n} \frac{|z|^2}{d^2}\frac{|f|^2}{d^{2(b+1)}}d\xi  \le \left(\int_{\mathbb{H}^n} \frac{|\nabla_{\mathbb{H}}f|^2}{d^{2(b+1)}}d\xi\right)^{\frac{1}{2}}
\left(\int_{\mathbb{H}^n} \frac{|z|^2}{d^2}\frac{|f|^2}{d^{2b}}d\xi\right)^{\frac{1}{2}}
\end{align*}
and $H^1_{b+1,b}(\mathbb{H}^n,\frac{|z|^2}{d^2}) \subset L^1_{b+1}(\mathbb{H}^n,\frac{|z|^2}{d^2})$.
\item If $f \in D^{1,2}({\mathbb{H}^n})$ then $f \in L^2_1(\mathbb{H}^n,\frac{|z|^2}{d^2})$ and
\begin{align*}
\left(\frac{Q-2}{2}\right)^2 \int_{\mathbb{H}^n} \frac{|z|^2}{d^2}\frac{|f|^2}{d^{2}} d\xi
\le \int_{\mathbb{H}^n} |\nabla_{\mathbb{H}}f|^2 d\xi.
\end{align*}
\item If $f \in H^1_{1,0}(\mathbb{H}^n,\frac{|z|^2}{d^2})$ then $f \in L^2_1(\mathbb{H}^n,\frac{|z|^2}{d^2})$ and 
\begin{align*}
\frac{Q-2}{2} \int_{\mathbb{H}^n} \frac{|z|^2}{d^2}\frac{|f|^2}{d^{2}} d\xi
\le \left(\int_{\mathbb{H}^n} \frac{|\nabla_{\mathbb{H}}f|^2}{d^2} d\xi\right)^{\frac{1}{2}} \left(\int_{\mathbb{H}^n} \frac{|z|^2}{d^2} |f|^2\xi\right)^{\frac{1}{2}}.
\end{align*}
\item If $f \in H^1_{a,-(a+1)}(\mathbb{H}^n,\frac{|z|^2}{d^2})$ then $f \in L^2(\mathbb{H}^n,\frac{|z|^2}{d^2})$ and 
\begin{align*}
\frac{Q}{2} \int_{\mathbb{H}^n} \frac{|z|^2}{d^2}|f|^2 d\xi
\le \left(\int_{\mathbb{H}^n} \frac{|\nabla_{\mathbb{H}}f|^2}{d^{2a}} d\xi\right)^{\frac{1}{2}} \left(\int_{\mathbb{H}^n} \frac{|z|^2}{d^2} |f|^2 d^{2(a+1)} d\xi\right)^{\frac{1}{2}}.
\end{align*}
\item If $f \in H^1_{-1,1}(\mathbb{H}^n,\frac{|z|^2}{d^2})$ then $f \in L^2_{\frac{1}{2}}(\mathbb{H}^n,\frac{|z|^2}{d^2})$ and 
\begin{align*}
\frac{Q-1}{2} \int_{\mathbb{H}^n} \frac{|z|^2}{d^2}\frac{|f|^2}{d} d\xi
\le \left(\int_{\mathbb{H}^n} |\nabla_{\mathbb{H}}f|^2 d^2 d\xi\right)^{\frac{1}{2}} \left(\int_{\mathbb{H}^n} \frac{|z|^2}{d^2} \frac{|f|^2}{d^2} d\xi\right)^{\frac{1}{2}}.
\end{align*}

\item If $f \in H^1(\mathbb{H}^n,\frac{|z|^2}{d^2})=H^1_{0,0}(\mathbb{H}^n,\frac{|z|^2}{d^2})$ then $f \in L^2_{\frac{1}{2}}(\mathbb{H}^n,\frac{|z|^2}{d^2})$ and 
\begin{align*}
\frac{Q-1}{2} \int_{\mathbb{H}^n} \frac{|z|^2}{d^2}\frac{|f|^2}{d} d\xi
\le \left(\int_{\mathbb{H}^n} |\nabla_{\mathbb{H}}f|^2 d\xi\right)^{\frac{1}{2}} \left(\int_{\mathbb{H}^n} \frac{|z|^2}{d^2} |f|^2 d\xi\right)^{\frac{1}{2}}.
\end{align*}

\item If $f \in H^1_{1,-1}(\mathbb{H}^n,\frac{|z|^2}{d^2})$ then $f \in L^2_{\frac{1}{2}}(\mathbb{H}^n,\frac{|z|^2}{d^2})$ and 
\begin{align*}
\frac{Q-1}{2} \int_{\mathbb{H}^n} \frac{|z|^2}{d^2} \frac{|f|^2}{d} d\xi
\le \left(\int_{\mathbb{H}^n} \frac{|\nabla_{\mathbb{H}}f|^2}{d^2} d\xi\right)^{\frac{1}{2}} \left(\int_{\mathbb{H}^n} |z|^2 |f|^2 d\xi\right)^{\frac{1}{2}}.
\end{align*}

\item If $f \in H^1_{-(b+1),b}(\mathbb{H}^n,\frac{|z|^2}{d^2})$ then $f \in L^2(\mathbb{H}^n,\frac{|z|^2}{d^2})$ and 
\begin{align*}
\frac{Q}{2} \int_{\mathbb{H}^n} \frac{|z|^2}{d^2}|f|^2 d\xi
\le \left(\int_{\mathbb{H}^n} |\nabla_{\mathbb{H}}f|^2 d^{2(b+1)} d\xi\right)^{\frac{1}{2}} \left(\int_{\mathbb{H}^n} \frac{|z|^2}{d^2} \frac{|f|^2}{d^{2b}} d\xi\right)^{\frac{1}{2}}.
\end{align*}
\end{enumerate}
\end{corollary}

\proof
By making special choices of $a,b \in \mathbb{R}$ in \eqref{e31} as follows:\\
(1) Let $b=a+1$  \ \ \ \ \ \ \ (4) Let $a=1$, $b+0$  \ \ \ \ \ \ (7) Let $a=0$, $b=0$\\
(2) Let $a=b+1$  \ \ \ \ \ \ \ (5) Let $b=-(a+1)$ \ \ \ \  \ (8) Let $a=1$, $b=-1$\\
(3) Let $a=0$, $b=1$ \ \ \ (6) Let $a=-1$, $b=1$ \ \ \ \ (9) Let $a=-(b+1)$. 

\qed

\subsection{Proof of Part (b) in Theorem \ref{thm32} and some consequences}   
In order to prove part (b) of Theorem \ref{thm32} we use the case $p=2$ in \eqref{e25} as stated in the next lemma.

\begin{lemma}\label{lem35}
For $Q\neq a+b-1$, $a,b \in \mathbb{R}$,
\begin{align}\label{e34}
\mathcal{D}_{a,b} \int_{\mathbb{H}^n} \frac{|z|^2}{d^2}\frac{|\nabla_{\mathbb{H}}f|^2}{d^{a+b+1}}d\xi  \le \left( \int_{\mathbb{H}^n} \frac{|\Delta_{\mathbb{H}}f|^2}{d^{2a}}d\xi\right)^{\frac{1}{2}}
\left(\int_{\mathbb{H}^n} \frac{|z|^2}{d^2}\frac{|\nabla_{\mathbb{H}}f|^2}{d^{2b}}d\xi\right)^{\frac{1}{2}}
\end{align}
with constant $\mathcal{D}_{a,b}:=\frac{|Q+a+b-1)|}{2}$.
\end{lemma}

We quickly remark that \eqref{e34} can also be proved  via Costa's approach \cite{Co2} as follows. For $s \in \mathbb{R}$ we have
\begin{align*}
\int_{\mathbb{H}^n} \left|\frac{\Delta_{\mathbb{H}}f}{d^a}+s\frac{\nabla_{\mathbb{H}}f}{d^b}\nabla_{\mathbb{H}^n} d \right|^2 d\xi \ge 0,
\end{align*}
that is,
\begin{align}\label{e35}
\int_{\mathbb{H}^n} \frac{|\Delta_{\mathbb{H}} f|^2}{d^{2a}}d\xi
+ s^2 \int_{\mathbb{H}^n} \frac{|z|^2}{d^2}\frac{|\nabla_{\mathbb{H}}f|^2}{d^{2b}}d\xi
+2s \Re{e}\int_{\mathbb{H}^n} \Delta_{\mathbb{H}} f\frac{\overline{\nabla_{\mathbb{H}}f} \nabla_{\mathbb{H}} d}{d^{a+b}}  d\xi \ge 0.
\end{align}
Applying the divergence theorem on the last integral in \eqref{e35}, we obtain
\begin{align}\label{e36a}
& \Re{e}\int_{\mathbb{H}^n} \text{div}_{\mathbb{H}} \nabla_{\mathbb{H}} f\frac{\overline{\nabla_{\mathbb{H}}f} \nabla_{\mathbb{H}} d}{d^{a+b}}  d\xi=  -\Re{e}  \int_{\mathbb{H}^n} \nabla_{\mathbb{H}} f\   \nabla_{\mathbb{H}}\cdot \left(\frac{\overline{\nabla_{\mathbb{H}}f} \nabla_{\mathbb{H}} d}{d^{a+b}} \right) d\xi \nonumber\\
& = -\frac{1}{2}  \int_{\mathbb{H}^n}  \nabla_{\mathbb{H}} ( |\nabla_{\mathbb{H}}f|^2)  \frac{ \nabla_{\mathbb{H}} d}{d^{a+b}}  d\xi - \int_{\mathbb{H}^n}   |\nabla_{\mathbb{H}}f|^2 \frac{ \nabla_{\mathbb{H}} (\nabla_{\mathbb{H}}d)}{d^{a+b}}  d\xi  + (a+b) \int_{\mathbb{H}^n}   |\nabla_{\mathbb{H}}f|^2 \frac{ |\nabla_{\mathbb{H}} d|^2}{d^{a+b+1}}  d\xi. 
\end{align}
The second term on the RHS of the last line also vanishes due to \eqref{e9a}, in the sense that 
$$ \nabla_{\mathbb{H}}(\nabla_{\mathbb{H}}d) = \frac{\frac{1}{2} \langle \nabla_\mathbb{H} |\nabla_\mathbb{H} d|^2, \nabla_\mathbb{H} d \rangle}{| \nabla_{\mathbb{H}} d|^2} = 0.$$
Applying the divergence theorem again, now on the first term on RHS of \eqref{e36a}, we have 
\begin{align}\label{e37a}
-\frac{1}{2} &  \int_{\mathbb{H}^n}  \nabla_{\mathbb{H}} ( |\nabla_{\mathbb{H}}f|^2)  \frac{ \nabla_{\mathbb{H}} d}{d^{a+b}}  d\xi  \nonumber\\
& =  \frac{1}{2}  \int_{\mathbb{H}^n}  |\nabla_{\mathbb{H}}f|^2  \frac{ \Delta_{\mathbb{H}} d}{d^{a+b}}  d\xi - \frac{a+b}{2}\int_{\mathbb{H}^n}   |\nabla_{\mathbb{H}}f|^2 \frac{ |\nabla_{\mathbb{H}} d|^2}{d^{a+b+1}}  d\xi.  
\end{align}
Combining \eqref{e36a} and \eqref{e37a} and using identities \eqref{e9} and \eqref{e10} we then arrive at 
\begin{align*}
\Re{e}\int_{\mathbb{H}^n} \text{div}_{\mathbb{H}} \nabla_{\mathbb{H}} f\frac{\overline{\nabla_{\mathbb{H}}f} \nabla_{\mathbb{H}} d}{d^{a+b}} d\xi= \frac{Q+a+b-1}{2} \int_{\mathbb{H}^n} \frac{|z|^2}{d^2}\frac{|\nabla_{\mathbb{H}}f|^2}{d^{a+b+1}}d\xi.
\end{align*}
Comparing \eqref{e35} with the inequality $As^2+Bs+C \ge 0$ and using the fact that $B^2-4AC\le 0$ holds true as before, then \eqref{cor34} follows.

\subsection*{Proof of Part (b) in Theorem \ref{thm32}}
As before Lemma \ref{lem35} implies
\begin{align*}
\int_{\mathbb{H}^n} \frac{|z|^2}{d^2}\frac{|\nabla_{\mathbb{H}}f|^2}{d^{a+b+1}}d\xi   \le \frac{2}{|Q-(a+b+1)|}\left(\int_{\mathbb{H}^n} \frac{|\Delta_{\mathbb{H}}f|^2}{d^{2a}}d\xi + \int_{\mathbb{H}^n} \frac{|z|^2}{d^2}\frac{|\nabla_{\mathbb{H}}f|^2}{d^{2b}}d\xi\right)
\end{align*}
for all $f \in C^\infty_0(\mathbb{H}^n\setminus\{0\})$. This implies the embedding b(i).  The embedding in b(ii)  follows from  \eqref{e34} by setting $b=a+1$, $a=\alpha \neq \frac{Q}{2}$ and $\alpha \le \frac{Q}{2}-2$, where $Q\ge a+b-1$.
\qed

As consequences of Lemma \ref{lem35} various forms of horizontal Hardy-Rellich and HPW type uncertainty principle inequalities can be  obtained.

\begin{corollary}\label{cor36}
The following inequalities hold true:
\begin{enumerate}[label=(\arabic*)]
\item For any $f \in H^2_{a,a+1}(\mathbb{H}^n,\frac{|z|^2}{d^2})$ and $a \le \frac{Q}{2}-2$,
\begin{align*}
\frac{|Q+2a|^2}{4} \int_{\mathbb{H}^n} \frac{|z|^2}{d^2}\frac{|\nabla_{\mathbb{H}}f|^2}{d^{2(a+1)}}d\xi  \le \int_{\mathbb{H}^n} \frac{|\Delta_{\mathbb{H}}f|^2}{d^{2a}}d\xi
\end{align*}
and $D^{2,2}_\alpha(\mathbb{H}^n) \subset D^{1,2}_{\alpha+1}(\mathbb{H}^n,\frac{|z|^2}{d^2})$.
\item For any $f \in D^{1,2}(\mathbb{H}^n,\frac{|z|^2}{d^2})$ and $Q\ge 3$, it holds that
\begin{align*}
\frac{Q-2}{2} \int_{\mathbb{H}^n} \frac{|z|^2}{d^2}\frac{|\nabla_{\mathbb{H}}f|^2}{d^{2(b+1)}}d\xi  \le \left(\int_{\mathbb{H}^n} \frac{|\Delta_{\mathbb{H}}f|^2}{d^{2a}}d\xi\right)^{\frac{1}{2}}
\left(\int_{\mathbb{H}^n} \frac{|z|^2}{d^2} |\nabla_{\mathbb{H}} f|^2d^{2(a+1)}d\xi\right)^{\frac{1}{2}}
\end{align*}
and $H^2_{a,-(a+1)}(\mathbb{H}^n,\frac{|z|^2}{d^2}) \subset D^{1,2}(\mathbb{H}^n,\frac{|z|^2}{d^2})$.

\item For $Q >2$ and any $f \in D^{1,2}_{\frac{1}{2}}(\mathbb{H}^n,\frac{|z|^2}{d^2})$ 
\begin{align*}
\frac{Q-1}{2} \int_{\mathbb{H}^n} \frac{|z|^2}{d^2}\frac{|\nabla_{\mathbb{H}}f|^2}{d} d\xi
\le \left(\int_{\mathbb{H}^n} \frac{|\Delta_{\mathbb{H}}f|^2}{d^{2a}} d\xi\right)^{\frac{1}{2}} \left(\int_{\mathbb{H}^n} |z|^2 |\nabla_{\mathbb{H}}f|^2 d\xi\right)^{\frac{1}{2}}.
\end{align*}

\item For $Q >2$ and any $f \in D^{1,2}_{\frac{1}{2}}(\mathbb{H}^n,\frac{|z|^2}{d^2})$ 
\begin{align*}
\frac{Q-1}{2} \int_{\mathbb{H}^n} \frac{|z|^2}{d^2} \frac{|f|^2}{d} d\xi
\le \left(\int_{\mathbb{H}^n} |\Delta_{\mathbb{H}}f|^2 d\xi\right)^{\frac{1}{2}} \left(\int_{\mathbb{H}^n} \frac{|z|^2}{d^2} |\nabla_{\mathbb{H}}f|^2 d\xi\right)^{\frac{1}{2}}.
\end{align*}

\item For any  $f \in D^{1,2}_1({\mathbb{H}^n},\frac{|z|^2}{d^2})$,
\begin{align*}
\left(\frac{Q}{2}\right)^2 \int_{\mathbb{H}^n} \frac{|z|^2}{d^2}\frac{|\nabla_{\mathbb{H}}f|^2}{d^{2}} d\xi
\le \int_{\mathbb{H}^n} |\Delta_{\mathbb{H}}f|^2 d\xi,
\end{align*}
$H^2_{0,1}(\mathbb{H}^n,\frac{|z|^2}{d^2}) \subset D^{1,2}_1({\mathbb{H}^n},\frac{|z|^2}{d^2})$ and $D^{2,2}_0({\mathbb{H}^n}) \subset D^{1,2}_1({\mathbb{H}^n},\frac{|z|^2}{d^2})$. 

\item For any  $f \in D^{1,2}_1({\mathbb{H}^n},\frac{|z|^2}{d^2})$,
\begin{align*}
\frac{Q}{2} \int_{\mathbb{H}^n} \frac{|z|^2}{d^2}\frac{\nabla_{\mathbb{H}}|f|^2}{d^{2}} d\xi
\le \left(\int_{\mathbb{H}^n} \frac{|\Delta_{\mathbb{H}}f|^2}{d^{2a}} d\xi\right)^{\frac{1}{2}}\left(\int_{\mathbb{H}^n} \frac{|z|^2}{d^2} |\nabla_{\mathbb{H}}f|^2\xi\right)^{\frac{1}{2}}.
\end{align*}
\end{enumerate}
\end{corollary}

\proof
By making special choices of $a,b \in \mathbb{R}$ in \eqref{e34} as follows:\\
(1) Let $b=a+1$  \ \ \ \ \ \ \ \ \ \ (3) Let $a=1$, $b=-1$  \  \ \ \ \ (5) Let $a=0$, $b=1$\\
(2) Let $b=-(a+1)$  \ \ \ \ \ \ (4) Let $a=0$, $b=0$ \ \ \ \ \   \ (6) Let $a=1$, $b=0$.
   
   \qed

\section{Inequalities for vector fields in domains of $\mathbb{H}^n$}\label{sec4}
 
In this section,  for an open subset $\Omega \subset \mathbb{H}^n$  we  let  $D^{1,p}(\Omega)$ denote the closure of $C^\infty_0(\Omega)$ with respect to the norm $\|f\|_{D^{1,p}(\Omega)} = (\int_\Omega |\nabla_{\mathbb{H}} f|^pd\xi)^{\frac{1}{p}}$.

 Let $\mathscr{W}^1(\Omega)$ be the space of locally integrable vector fields $V$ on $\Omega$ for which $V$ has weak horizontal divergence. By the way of definition, the function $A_V \in L^1_{loc}(\Omega)$ is  called the weak horizontal divergence of $V$ if and only if $\text{div}_{\mathbb{H}}V=A_V$ in the sense of distribution, that is, 
\begin{align}\label{e41}
\int_\Omega \phi A_V d\xi = - \int_\Omega \nabla_{\mathbb{H}} \phi \cdot V d\xi
\end{align}
for every $\phi \in C^\infty_0(\Omega)$. If we pick any vector field $V \in  \mathscr{W}^1(\Omega)$ and $f \in C^\infty_0(\Omega)$, then the vector field $|f|^pV$ has compact support in $\Omega$ and 
\begin{align}\label{e42}
\int_\Omega \text{div}_{\mathbb{H}} (|f|^pV)d\xi =0.
\end{align}

\begin{proposition}\label{pro41}
Let $\Omega$ be a bounded open subset of $\mathbb{H}^n$.   For  $p>1$ and any   vector field $V \in \mathscr{W}^1(\Omega)$,
\begin{align}\label{e43}
\int_\Omega |\nabla_{\mathbb{H}} f|^pd\xi \ge \int_\Omega\left[\text{div}_{\mathbb{H}} V - (p-1)|V|^{\frac{p}{p-1}} \right]|f|^p d\xi.
\end{align}
\end{proposition}

This result can be extended to the whole of $\mathbb{H}^n$ as a result of the embedding $D^{1,p}(\mathbb{H}^n) \subset D^{1,p}(\Omega)$.
\begin{proposition}
For any vector field $V \in \mathscr{W}^1(\mathbb{H}^n)$ and $p>1$,
\begin{align}\label{e44}
\int_{\mathbb{H}^n} |\nabla_{\mathbb{H}} f|^pd\xi \ge \int_{\mathbb{H}^n}\left[\text{div}_{\mathbb{H}} V - (p-1)|V|^{\frac{p}{p-1}} \right]|f|^p d\xi.
\end{align}
\end{proposition}

\proof
For any vector field $V \in \mathscr{W}^1(\Omega)$ and $f \in D^{1,p}(\Omega)$, it follows that $|f|^pV \in \mathscr{W}^1$ and we have 
\begin{align}\label{e45}
\text{div}_{\mathbb{H}} (|f|^pV) = \nabla_{\mathbb{H}} |f|^p V + |f|^p \text{div}_{\mathbb{H}} V.
\end{align}
Combining \eqref{e45} with the divergence theorem \eqref{e42} we get
\begin{align}\label{e46}
\int_{\Omega} |f|^p \text{div}_{\mathbb{H}} Vd\xi &= - p \int_{\Omega} |f|^{p-1}\nabla_{\mathbb{H}}|f|V d\xi \nonumber\\
& \le p \int_{\Omega} |f|^{p-1} |\nabla_{\mathbb{H}}f||V| d\xi\\
&\le p\left( \int_{\Omega} |\nabla_{\mathbb{H}}f|^pd\xi \right)^{\frac{1}{p}} \left(\int_{\Omega}|V|^{\frac{p}{p-1}}|f|^pd\xi\right)^{\frac{p-1}{p}},\nonumber
\end{align}
where we used the fact that $|\nabla_{\mathbb{H}}|f||\le |\nabla_{\mathbb{H}}f|$ and the H\"older inequality.

Consideration of Young's inequality of the form $\Phi\Psi \le \frac{1}{p}\Phi^p+\frac{(p-1)}{p}\Psi^{\frac{p}{p-1}}$ and setting $\Phi:=\left( \int_{\Omega} |\nabla_{\mathbb{H}}f|^pd\xi \right)^{\frac{1}{p}}$ and $\Psi:= \left(\int_{\Omega}|V|^{\frac{p}{p-1}}|f|^pd\xi\right)^{\frac{p-1}{p}}$ leads to 
\begin{align}\label{e47}
p\left( \int_{\Omega} |\nabla_{\mathbb{H}}f|^pd\xi \right)^{\frac{1}{p}} &  \left(\int_{\Omega}|V|^{\frac{p}{p-1}}|f|^pd\xi\right)^{\frac{p-1}{p}} \nonumber\\
& \le \int_{\Omega} |\nabla_{\mathbb{H}}f|^pd\xi +(p-1)\int_{\Omega}|V|^{\frac{p}{p-1}}|f|^pd\xi.
\end{align}
Putting \eqref{e46} and \eqref{e47} together and rearranging gives 
\begin{align*}
\int_\Omega |\nabla_{\mathbb{H}} f|^pd\xi \ge \int_\Omega\left[\text{div}_{\mathbb{H}} V - (p-1)|V|^{\frac{p}{p-1}} \right]|f|^p d\xi.
\end{align*}
This completes the proof.
\qed

In the next part we specialize $V$ into  various quantities and obtain specific weighted Hardy inequalities.

\begin{example}
Choose $$V=\alpha\frac{|\nabla_{\mathbb{H}}d|^{p-2}\nabla_{\mathbb{H}}d}{d^{p-1}},$$
where $\alpha$ is to be determined.
Clearly $$|V|^{\frac{p}{p-1}} = |\alpha|^{\frac{p}{p-1}}\frac{|\nabla_{\mathbb{H}}d|^p}{d^p} = |\alpha|^{\frac{p}{p-1}} \frac{|z|^p}{d^{2p}}$$
and
$$\text{div}_{\mathbb{H}} V = \alpha \frac{\nabla_{\mathbb{H}}(|\nabla_{\mathbb{H}}d|^{p-2})\nabla_{\mathbb{H}} d}{d^{p-1}}
+\alpha\frac{|\nabla_{\mathbb{H}}d|^{p-2}\Delta_{\mathbb{H}} d}{d^{p-1}} - \alpha(p-1)\frac{|\nabla_{\mathbb{H}}d|^p}{d^p}.$$
Using \eqref{e10}, \eqref{e11} and \eqref{e12} we obtain for $p\neq Q$ that
$$\text{div}_{\mathbb{H}} V  = - \alpha(p-Q)\frac{|z|^p}{d^{2p}}.$$
Thus
$$\text{div}_{\mathbb{H}} V - (p-1)|V|^{\frac{p}{p-1}} = -(p-1)\left[\alpha\frac{p-Q}{p-1}+|\alpha|^{\frac{p}{p-1}}\right]\frac{|z|^p}{d^{2p}}.$$
Therefore \eqref{e44} yields
\begin{align}\label{e48}
\int_{\mathbb{H}^n} |\nabla_{\mathbb{H}} f|^pd\xi \ge - (p-1)\left[\alpha\frac{p-Q}{p-1}+|\alpha|^{\frac{p}{p-1}}\right] \int_{\mathbb{H}^n} \frac{|z|^p}{d^{2p}}|f|^p d\xi.
\end{align}
Optimizing the quantity $\alpha\mapsto F(\alpha):=-(p-1)\left[\alpha\frac{p-Q}{p-1}+|\alpha|^{\frac{p}{p-1}}\right]$, we obtain $\frac{p-Q}{p-1}+\frac{p}{p-1}|\alpha|^{\frac{1}{p-1}}=0$ giving the extremal value $\alpha_\star =-\left(\frac{p-Q}{p}\right)^{p-1}$ and then
$$F(\alpha_\star)=\left(\frac{Q-p}{p}\right)^p.$$
By this, \eqref{e48} yields the weighted inequality
\begin{align*}
\int_{\mathbb{H}^n} |\nabla_{\mathbb{H}} f|^pd\xi \ge \left(\frac{Q-p}{p}\right)^p \int_{\mathbb{H}^n} \frac{|z|^p}{d^p}\frac{|f|^p}{d^p} d\xi.
\end{align*}
The above is another approach to proving \eqref{e23}. Generalised Picone identity has also been applied by Niu, Zhang and Wang \cite{NZW} to derive this inequality.
\end{example}

\begin{example}
Let $\mathbb{H}^+:=\{\xi \in \mathbb{H}^n : \langle \xi,\nu\rangle >\delta\}$ be the half space of the Heisenberg group $\mathbb{H}^n$, where $\nu =(\nu_x,\nu_y,\nu_t)$  with $\nu_x,\nu_y \in \mathbb{R}^n$ and $\nu_t\in \mathbb{R}$ becomes the Riemannian outer unit normal to the boundary $\partial \mathbb{H}^+$ and $\delta\in\mathbb{R}$ (\cite{Gar,RSS}).   Then for all $f \in C^\infty_0(\mathbb{H}^+)$ and $p>1$, we have 
\begin{align}\label{e49}
\int_{\mathbb{H}^+} |\nabla_{\mathbb{H}} f|^pd\xi \ge \left(\frac{p-1}{p}\right)^p \int_{\mathbb{H}^+} \frac{W(\xi)^p}{\text{dist}(\xi,\partial \mathbb{H}^+)^p} |f|^p d\xi
\end{align} 
where $W(\xi):=(\sum_{j=1}^n\langle X_j(\xi),\nu\rangle^2+ \langle Y_j(\xi),\nu\rangle^2)^{1/2}$ is the so-called angle function,
and $\text{dist}(\xi,\partial \mathbb{H}^+) := \langle \xi,\nu\rangle-d$ is the Euclidean distance to the boundary $\partial \mathbb{H}^+$.

To obtain \eqref{e49} we can apply \eqref{e43} by choosing $V$  as 
$$V = \beta\frac{|\nabla_{\mathbb{H}} \widetilde{\delta}|^{p-2}\nabla_{\mathbb{H}} \widetilde{\delta}}{\widetilde{\delta}^{p-1}}$$
with $\beta:=-(\frac{p-1}{p})^{p-1}$, where $\widetilde{\delta}$ defined as $\widetilde{\delta}:=\text{dist}(\xi,\partial \mathbb{H}^+)$.

Note that inequality \eqref{e49} was conjectured by Larson \cite{Lar} and later proved by Ruzhansky, Sabitbek and Suragan \cite[Corollary 2.2]{RSS}. Also, the sharpness of the inequality  was proposed by Larson \cite{Lar} by choosing $\nu := (1,0,\cdots,0)$ and $\delta = 0$. 
\end{example}

\begin{example}
Let $\Omega =\mathscr{B}_R\setminus\{0\}$, where $\mathscr{B}_R:=\{\xi \in \mathbb{H}^n : d(\xi,0)<R\}$ is the Heisenberg ball, and let  $\widetilde{\delta}:=R-d>0.$
Choose
$$V=  - \alpha\frac{|\nabla_{\mathbb{H}}\widetilde{\delta}|^{p-2} \nabla_{\mathbb{H}}\widetilde{\delta}}{\widetilde{\delta}^{p-1}}$$
and $\alpha = (\frac{p-1}{p})^{p-1}$.
One then computes
$$|V|^{\frac{p}{p-1}}= |\alpha|^{\frac{p}{p-1}} \frac{|\nabla_\mathbb{H}d|^p}{\widetilde{\delta}^p}=  \left(\frac{p-1}{p}\right)^p\frac{|z|^p}{d^p \widetilde{\delta}^p}$$
and
\begin{align*}
\text{div}_{\mathbb{H}} V &= \alpha \frac{\nabla_{\mathbb{H}}(|\nabla_{\mathbb{H}}d|^{p-2}\nabla_{\mathbb{H}} d)}{\widetilde{\delta}^{p-1}}
 + \alpha(p-1)\frac{|\nabla_{\mathbb{H}}d|^p}{\widetilde{\delta}^p} \\
 &= \alpha\frac{|\nabla_{\mathbb{H}}d|^{p-2}\Delta_{\mathbb{H}} d}{\widetilde{\delta}^{p-1}} + \alpha(p-1)\frac{|\nabla_{\mathbb{H}}d|^p}{\widetilde{\delta}^p}\\
 & = \alpha\left[(Q-1)\frac{\delta}{d}+(p-1) \right]\frac{|z|^p}{d^p \widetilde{\delta}^p} 
 \ge \alpha (p-1)\frac{|z|^p}{d^p \widetilde{\delta}^p},
\end{align*}
where we have used identities \eqref{e9}, \eqref{e10}  and \eqref{e12}. 
Then applying \eqref{e43} we have  for $p>1$
\begin{align}\label{e410}
\int_{\Omega} |\nabla_{\mathbb{H}} f|^pd\xi \ge \left(\frac{p-1}{p}\right)^p \int_{\Omega} \frac{|z|^p}{d^p}\frac{|f|^p}{\widetilde{\delta}^p} d\xi.
\end{align}
Inequality \eqref{e410}  has been derived by Huan and Niu \cite[Theorem 2.3]{HN} through another method.
\end{example}

\begin{example}
Let $\Omega \subseteq \mathbb{H}^n$ be bounded.
Choose
$$V= \gamma^{p-1} \frac{|\nabla_{\mathbb{H}}d|^{p-2} \nabla_{\mathbb{H}} d}{(d\ln R/d)^{p-1}},$$
where $0<\sup_\Omega d<R<+\infty$ and $\gamma = \frac{p-1}{p}$.
If $1<p=Q$
\begin{align*}
\text{div}_\mathbb{H}V &= \gamma^{p-1}\left[\left(\frac{\Delta_{\mathbb{H},p}d}{d^{p-1}}-(p-1)\frac{|\nabla_\mathbb{H}d|^p}{d^p}\right)(\ln R/d)^{1-p} +(p-1)\frac{|\nabla_\mathbb{H}d|^p}{d^p}(\ln R/d)^{-p} \right]\\
& = \gamma^{p-1}\left[\left(\frac{Q-1-(p-1)}{d^p}\right)\frac{|z|^p}{d^p}(\ln R/d)^{1-p}  +(p-1)\frac{|\nabla_\mathbb{H}d|^p}{d^p}(\ln R/d)^{-p} \right]\\
& = (p-1)\gamma^{p-1} \frac{|\nabla_\mathbb{H}d|^p}{d^p(\ln R/d)^p}.
\end{align*}
Therefore
\begin{align*}
\text{div}_\mathbb{H} V -(p-1)|V|^{\frac{p}{p-1}} &= (p-1)\left[\gamma^{p-1} -|\gamma|^p \right] \frac{|\nabla_\mathbb{H}d|^p}{d^p(\ln R/d)^p} \\
& =  \left(\frac{p-1}{p}\right)^p \frac{|z|^p}{d^p}\frac{1}{(d\ln R/d)^p}.
\end{align*}

Hence, by \eqref{e43} we obtain the following improved weighted Hardy inequality for $f \in C^\infty_0(\Omega)$ and $1<p=Q$
\begin{align*}
\int_{\Omega} |\nabla_{\mathbb{H}} f|^pd\xi \ge \left(\frac{p-1}{p}\right)^p \int_{\Omega} \frac{|z|^p}{d^p}\frac{|f|^p}{d^p(\ln R/d)^p} d\xi.
\end{align*}
One can see paper \cite{DNY} by Dou, Niu and Yuan for variants of this.
\end{example}

 
\subsection*{Acknowledgement} 
This paper was completed during the first author's research visit to Ghent Analysis and PDE Centre, Ghent University.  He therefore gratefully acknowledges the  research supports of  IMU-Simons African Fellowship Grant  and EMS-Simons for African  program.  He also thanks his host Professor Michael Ruzhansky.  This research has been funded by the FWO Odysseus 1 grant G.0H94.18N: Analysis and Partial Differential Equations and by the Methusalem programme of the Ghent University Special Research Fund (BOF) (Grant number 01M01021).  Michael Ruzhansky is also supported by EPSRC grant EP/R003025/2. 


\end{document}